\documentclass[12pt]{article}

\usepackage{amsmath,amssymb,amsfonts,graphicx}
\usepackage[latin1]{inputenc}

\newtheorem{thm}{Theorem}[section]
\newtheorem{prop}[thm]{Proposition}

\newtheorem{df}[thm]{Definition}

\newtheorem{rem}[thm]{Remark}
\newtheorem{ass}[thm]{Assumption}
\newtheorem{cor}[thm]{Corollary}

\def\be#1 {\begin{equation} \label{#1}}
\newcommand{\ee}{\end{equation}}
\def\dem {\noindent {\bf Proof : }}

\newcommand{\mb}{\medskip\noindent}
\newcommand{\gb}{\bigskip\noindent}
\newcommand{\R}{\mathbb R}

\newcommand{\p}{s}

\newcommand{\D}{\mathcal D}

\newcommand{\diam}{diam}
\newcommand{\aver}[1]{-\hskip-0.46cm\int_{#1}}

\newcommand{\ind}{ {\bf 1} }

\DeclareMathOperator{\Lip}{Lip}

\def\sqw{\hbox{\rlap{\leavevmode\raise.3ex\hbox{$\sqcap$}}$ \sqcup$}}
\def\findem{\ifmmode\sqw\else{\ifhmode\unskip\fi\nobreak\hfil
\penalty50\hskip1em\null\nobreak\hfil\sqw
\parfillskip=0pt\finalhyphendemerits=0\endgraf}\fi}

\title{Maximal inequalities for dual Sobolev spaces $W^{-1,p}$ and applications to interpolation.}

\author{Fr\'ed\'eric Bernicot \\ frederic.bernicot@math.u-psud.fr \\ D\'epartement de Math\'ematiques \\ Universit\'e Paris-Sud \\ 91405 Orsay Cedex (France)}

\begin {document}

\maketitle

\begin{abstract} We firstly describe a maximal inequality for dual Sobolev spaces $W^{-1,p}$. This one corresponds to a ``Sobolev version'' of usual properties of the Hardy-Littlewood maximal operator in Lebesgue spaces. Even in the euclidean space, this one seems to be new and we develop arguments in the general framework of Riemannian manifold. Then we present an application to obtain interpolation results for Sobolev spaces.
\end{abstract}

\mb {\bf Key-words~:} Maximal inequalities, Sobolev spaces, interpolation. 

\mb {\bf MSC~:} 46E35-42B25-46B70.

\tableofcontents

\gb
The first maximal inequality in Lebesgue spaces, is described by the $L^p$-boundedness of the Hardy-Littlewood maximal function. This result holds in a space of homogeneous type $(X,d,\mu)$ ~: for $p\in(1,\infty]$, $s\in[1,p)$ and $f\in L^p(X)$
$$ \|f\|_{L^p(X)} \lesssim \left\| x\to \sup_{ \genfrac{}{}{0pt}{}{Q\ \textrm{ball}}{Q\ni x} } \frac{1}{\mu(Q)^{1/s}} \left\| f\right\|_{L^s(Q)} \right\|_{L^p(X)} \lesssim \|f\|_{L^p(X)}.$$
Here the left inequality is due to the ``regularity property''~: for almost every $x\in X$
\be{pro} \lim_{r\to 0} \frac{1}{\mu(B(x,r))} \int_{B(x,r)} |f| d\mu = |f(x)|. \ee
The right one corresponds to the $L^p$ boundedness of the maximal operator.

\mb Applying this result to a function and its gradient, we obtain the same result for the Sobolev spaces on a doubling Riemannian manifold $M$~: for $p\in(1,\infty]$, $s\in[1,p)$ and $f\in W^{1,p}$
$$ \|f\|_{W^{1,p}} \lesssim \left\| x\to \sup_{ \genfrac{}{}{0pt}{}{Q\ \textrm{ball}}{Q\ni x} } \frac{1}{\mu(Q)^{1/s}} \left\| f\right\|_{W^{1,s}(Q)} \right\|_{L^p} \lesssim \|f\|_{W^{1,p}}.$$
Therefore the Sobolev norm can easily be described by the corresponding Lebesgue norm of a maximal operator (which is a ``Sobolev version'' of the Hardy-Littlewood maximal function). Such a property is important because the norms in Lebesgue spaces are specific and satisfies for example the ``lattice property'' which is not the case of the norms in Sobolev spaces.
 Then a natural question arises~: do we have similar results for the dual Sobolev spaces $W^{-1,p}$ ?

\gb Recently in \cite{BJ,B2} the authors have used maximal operators (and duality) to describe interpolation results between Hardy spaces and Lebesgue spaces. To extend this theory for Sobolev spaces, we need such maximal inequalities for negative Sobolev spaces.
That is why, we study this problem. Despite this objective, the above inequality studied in the current paper may be of independent interest by itself.

\gb We define maximal operators and then prove the following result~: under classical assumptions on the Riemannian manifold $M$, there are implicit constants such that for all functions $f\in W^{-1,p}$ 
\be{eq:maxx} \|f\|_{W^{-1,p}} \underset{(a)}{\lesssim} \left\| \sup_{ \genfrac{}{}{0pt}{}{Q\ \textrm{ball}}{x\in Q} } \frac{1}{\mu(Q)^{1/s}} \left\| f\right\|_{W^{-1,s}(Q)} \right\|_{L^p} \underset{(b)}{\lesssim} \|f\|_{W^{-1,p}}\ee
under some restrictions on $p$.
The second inequality $(b)$ is quite easy to obtain and corresponds to a boundedness of the maximal operator. The first one $(a)$ is more difficult to prove. Such property as (\ref{pro}) is not sufficient to conclude. \\
For example in the euclidean space, we get~:

\begin{thm} On ${\mathbb R}^n$ equipped with the euclidean structure, (\ref{eq:maxx}) holds for every exponents $p,r\in(1,\infty)$.
\end{thm}

We emphasize that even in the euclidean space $\R^n$, such inequalities are not obvious. In this particular case, we know that the operator $(I+\Delta)^{-1/2}$ defines an isomorphism from $L^p(\R^n)$ to $W^{-1,p}(\R^n)$. However such a description is not sufficiently precise to obtain the inequality $(a)$. \\
This result seems to be new and does not exist in the litterature. We think that it will permit to better understand the structure of dual Sobolev spaces and above all the interactions with restriction and localization operators. \\
We believe in the interest of such inequalities and we give a first application about interpolation of Sobolev spaces (Section \ref{sec:interpolation}). For example we will prove the following result.

\begin{thm} Let $M$ be a doubling Riemannian manifold satisfying a {\em Reverse Riesz inequality}~:
$$ \left\|(1+\Delta)^{1/2}(f) \right\|_{L^r} \lesssim \|f\|_{W^{1,r}},$$ 
for an exponent $r\in(1,2)$. 
Then for all $p_0\in(1,2)$ and $\theta\in(0,1)$ such that
$$ \frac{1}{p_\theta}:=\frac{\theta}{p_0}+\frac{1-\theta}{2} <\frac{1}{r},$$
we have
$$\left(W^{1,p_0},W^{1,2}\right)_{\theta,p_\theta}=W^{1,p_\theta}.$$
\end{thm}

\mb This result is interesting as we do not require Poincar\'e inequality as in the work of N. Badr (see \cite{Nadine,Nadine1}). This is the first result of interpolation for Sobolev spaces, which permits to get around the use of Poincar\'e inequalities. Due to the work of P.~Auscher and T.~Coulhon (see \cite{AC}), our assumed {\em Reverse Riesz inequality} is weaker than the Poincar\'e inequality $(P_{r})$.

\gb We refer the reader to a forthcoming work (joined with N. Badr, see \cite{HSNadine}), where we use these maximal inequalities for Sobolev spaces in order to describe an interpolation theory for abstract Hardy-Sobolev spaces. In this case, they will play a crucial role.

\section{Preliminaries}

Throughout this paper we will denote by $\ind_{E}$ the characteristic function of
 a set $E$ and $E^{c}$ the complement of $E$. If $X$ is a metric space, $\Lip$ will be  the set of real Lipschitz functions on $X$ and $\Lip_{0}$ the set of real, compactly supported Lipschitz functions on $X$. For a ball $Q$ in a metric space, $\lambda Q$  denotes the ball co-centered with $Q$ and with radius $\lambda$ times that of $B$. Finally, $C$ will be a constant that may change from an inequality to another and we will use $u\lesssim
v$ to say that there exists two constants $C$  such that $u\leq Cv$ and $u\simeq v$ to say that $u\lesssim v$ and $v\lesssim u$.

\mb In all this paper $M$ denotes a Riemannian manifold. We write $\mu$ for the Riemannian measure on $M$, $\nabla$ for the
Riemannian gradient, $|\cdot|$ for the length on the tangent space (forgetting the subscript $x$ for simplicity) and
$\|\cdot\|_{L^p}$ for the norm on $ L^p:=L^{p}(M,\mu)$, $1 \leq p\leq +\infty.$  We denote by $Q(x, r)$ the open ball of
center $x\in M $ and radius $r>0$. \\
We will use the positive Laplace-Beltrami operator $\Delta$ defined by
$$ \forall f,g\in C^\infty_0(M), \qquad \langle \Delta f,g\rangle = \langle \nabla f ,\nabla g \rangle.$$

\subsection{The doubling property}

\begin{df} Let $M$ be a Riemannian manifold. One says that $M$ is {\em doubling} or satisfies the (global) doubling property $(D)$ if there exists a
constant $C>0$, such that for all $x\in M,\, r>0 $ we have
\begin{equation*}\tag{$D$}
\mu(Q(x,2r))\leq C \mu(Q(x,r)).
\end{equation*}
\end{df}
\noindent Observe that if $M$ satisfies $(D)$ then
$$ \diam(M)<\infty\Leftrightarrow\,\mu(M)<\infty\,\textrm{ (see \cite{ambrosio1})}. $$
Therefore if $M$ is a complete Riemannian manifold satisfying $(D)$ then $\mu(M)=\infty$.

\begin{thm}[Maximal theorem]\label{MIT} (\cite{coifman2})
Let $M$ be a Riemannian manifold satisfying $(D)$. Denote by $M_{HL}$ the uncentered Hardy-Littlewood maximal function
over open balls of $M$ defined by
 $$ M_{HL}f(x):=\underset{\genfrac{}{}{0pt}{}{Q \ \textrm{ball}}{x\in Q}} {\sup} \ \frac{1}{\mu(Q)}\int_{Q}|f| d\mu.$$
Then for every  $p\in(1,\infty]$, $M_{HL}$ is $L^p$ bounded and moreover of weak type $(1,1)$.
%\footnote{ An operator $T$ is of weak type $(p,p)$ if there is $C>0$ such that for any $\alpha>0$, $\mu(\{x;\,|Tf(x)|>\alpha\})\leq \frac{C}{\alpha^p}\|f\|_p^p$.}.
\\
Consequently for $s\in(0,\infty)$, the operator $M_{HL,s}$ defined by
$$ M_{HL,s}f(x):=\left[M_{HL}(|f|^s)(x) \right]^{1/s} $$
is of weak type $(s,s)$ and $L^p$ bounded for all $p\in(s,\infty]$.
\end{thm}

\subsection{Poincar\'e inequality}
\begin{df}[Poincar\'{e} inequality on $M$] We say that a Riemannian manifold $M$ admits \textbf{a Poincar\'{e} inequality $(P_{q})$} for some $q\in[1,\infty)$ if there exists a constant $C>0$ such that, for every function $f\in \Lip_{0}(M)$ and every ball $Q$ of $M$ of radius $r>0$, we have
\begin{equation*}\tag{$P_{q}$}
\left(\aver{Q}|f-f_{Q}|^{q} d\mu\right)^{1/q} \leq C r \left(\aver{Q}|\nabla f|^{q}d\mu\right)^{1/q}.
\end{equation*}
\end{df}
\begin{rem} By density of $C_{0}^{\infty}(M)$ in $Lip_0(M)$, we can replace $\Lip_{0}(M)$ by $C_{0}^{\infty}(M)$.
\end{rem}

\mb Let us recall some known facts about Poincar\'{e} inequalities with varying $q$.
 \\
It is known that $(P_{q})$ implies $(P_{p})$ when $p\geq q$ (see \cite{hajlasz4}). Thus if the set of $q$ such that
$(P_{q})$ holds is not empty, then it is an interval unbounded on the right. A recent result of S. Keith and X. Zhong
(see \cite{KZ}) asserts that this interval is open in $[1,+\infty[$~:

\begin{thm}\label{kz} Let $(X,d,\mu)$ be a complete metric-measure space with $\mu$ doubling
and admitting a Poincar\'{e} inequality $(P_{q})$, for  some $1< q<\infty$.
Then there exists $\epsilon >0$ such that $(X,d,\mu)$ admits
$(P_{p})$ for every $p>q-\epsilon$.
\end{thm}

%\mb A consequence of Poincar\'e inequality:

\section{Maximal characterization of dual Sobolev spaces} \label{sec:maximal}

From now on, we always assume that the Riemannian manifold satisfies the doubling property $(D)$ and write $n$ for its dimension.

\subsection{New maximal operators.}

\mb First, we begin recalling the ``duality-properties'' of the Sobolev spaces.

\begin{df} For $p\in[1,\infty]$ and $O$ an open set of $M$, we define $W^{1,p}(O)$ as following
$$ W^{1,p}(O) := \overline{C^\infty_0(O)}^{\|\ \|_{W^{1,p}(O)}} \quad \textrm{with} \quad \|f\|_{W^{1,p}(O)}:=\left\||f|+|\nabla f| \right\|_{L^p(O)}.$$
Then we denote $W^{-1,p'}(O)$ the dual space of $W^{1,p}(O)$ defined as the set of distributions $f\in \mathcal{D}'(M)$ such that
$$
\|f\|_{W^{-1,p'}(O)}=\sup_{g\in C^{\infty}_{0}(M)} \frac{|\langle f,g\rangle|}{\|g\|_{W^{1,p}(O)}}.
$$
\end{df}

\begin{prop} \label{dualS} Let $p\in[1,\infty)$. Then for all open set $O$ of $M$, we have
\begin{align*}
 \|f\|_{W^{-1,p'}(O)} & = \inf_{f=\phi-div(\psi)} \left\|\phi\right\|_{L^{p'}(O)} + \left\| \psi\right\|_{L^{p'}(O)} \\
 & \simeq \inf_{f=\phi-div(\psi)} \left\||\phi|+|\psi|\right\|_{L^{p'}(O)}.
\end{align*}
Here we take the infimum over all the decompositions $f=\phi-div(\psi)$ on $M$ with $\phi\in L^{p'}(O)$ and $\psi\in
\D'(O,\R^n)$ such that $div(\psi)\in L^{p'}(O)$.
\end{prop}

\mb The proof is left to the reader (it is essentially written in \cite{ART}, Proposition 33).

\mb We introduce the following maximal operators~:

\begin{df} Let $s>0$. According to the standard maximal ``Hardy-Littlewood'' operator
$ M_{HL,s}$,
we define two ``Sobolev versions''~:
$$ M_{S,s}(f)(x):=\sup_{\genfrac{}{}{0pt}{}{Q \textrm{ball}}{Q\ni x}} \ \frac{1}{\mu(Q)^{1/s}} \left\|f \right\|_{W^{-1,s}(Q)}$$
and
$$ M_{S,*,s}(f)(x):= \inf_{f = \phi - div(\psi)} \ M_{HL,s}\left(|\phi|+|\psi|\right)(x).$$
\end{df}

\begin{rem} \label{rem:inegalite} Thanks to Proposition \ref{dualS}, it is easy to check that we can compare them pointwisely~:
$$ M_{S,s}(f) \leq M_{S,*,s}(f).$$
\end{rem}

\mb We dedicate the next subsection to the study of these maximal operators. Mainly we want to describe the dual Sobolev norms by the corresponding Lebesgue norms of these operators.

\subsection{First properties of the maximal operators.} \label{subsec:maximal}

\mb We begin proving some useful and general properties for the new maximal operators $M_{S,s}$ and $M_{S,*,s}$. These
operators can be thought as being equivalent to $M_{HL,s}((I+\Delta)^{-1/2})$, where $\Delta$ is the positive
Laplace-Beltrami operator on the manifold $M$.

\begin{prop} \label{prop:importante}
For $p\in[1,\infty)$, $M_{S,p}$ and $M_{S,*,p}$ are of ``weak type $(p,p)$''~: there exists an implicit constant such that for all $f\in W^{-1,p}$ \be{eq1} \left\| M_{S,p}(f) \right\|_{L^{p,\infty}} \leq \left\| M_{S,*,p}(f) \right\|_{L^{p,\infty}} \lesssim
\|f\|_{W^{-1,p}}. \ee
\end{prop}

\dem The first inequality is due to Remark \ref{rem:inegalite}. We only check the second one. Using Fatou's lemma in weak Lebesgue spaces, it yields
$$\left\| M_{S,*,p}(f) \right\|_{L^{p,\infty}} \leq \inf_{f=\phi-div(\psi)} \left\| M_{HL,p}(|\phi|+|\psi|) \right\|_{L^{p,\infty}}.$$
Then using the weak type $(p,p)$ of the Hardy-Littlewood maximal operator it comes
$$\left\| M_{S,*,p}(f) \right\|_{L^{p,\infty}} \lesssim \inf_{f=\phi-div(\psi)} \left\||\phi|+|\psi|\right\|_{L^p}.$$
Finally Proposition  \ref{dualS} finishes the proof. \findem

\mb Now we look for reverse inequalities. First we describe an easy fact~:

\begin{rem} \label{rem:comparaison} Let $r_1,r_2\in [1,\infty)$ with $r_1\leq r_2$. Then
 $$M_{S,*,r_1}\leq M_{S,*,r_2} \quad \textrm{and} \quad M_{S,r_1}\leq M_{S,r_2}.$$
\end{rem}

\begin{prop} \label{prop:importante2} Let $p\in [1,\infty)$. The two maximal operators $M_{S,*,p}$ and $M_{S,p}$ ``control the Sobolev norm in $W^{-1,p}$''. That is
\be{eq3bis} \forall f\in W^{-1,p}, \qquad \|f\|_{W^{-1,p}} \lesssim \left\| M_{S,p}(f) \right \|_{L^p} \leq \left\| M_{S,*,p}(f) \right \|_{L^p}. \ee
\end{prop}

\dem Thanks to Remark \ref{rem:inegalite}, we just have to prove the first inequality. In order to show this one, we choose a collection of balls $(B_i)_i$ of radius $1$, which corresponds to a bounded covering of $M$. Let $(\phi_i)_i$ be a partition of unity associated to this covering. Then we know that there exists a function $g\in C^\infty_0$ such that
$$ \| f\|_{W^{-1,p}} \leq 2 \langle f, g \rangle = 2\sum_i \langle f, g\phi_i \rangle$$
and $\|g\|_{W^{1,p'}}=1$. We use the fact that
$$ \langle f, g\phi_i \rangle \leq \|f\|_{W^{-1,p}(Q_i)} \|g\phi_i\|_{W^{1,p'}(Q_i)}.$$
Since the balls $B_i$ are of radius $1$, the functions $\phi_i$ can be chosen as uniformly bounded in the Sobolev space $W^{1,p'}$ and so we have
$$ \langle f, g\phi_i \rangle \lesssim \|f\|_{W^{-1,p}(Q_i)} \|g\|_{W^{1,p'}(Q_i)} \lesssim \mu(Q_i)^{1/p} \inf_{Q_i} M_{S,p}(f) \|g\|_{W^{1,p'}(Q_i)}.$$
Using H\"older inequality we obtain
$$ \| f\|_{W^{-1,p}} \lesssim  \left(\sum_{i} \mu(Q_i) \inf_{Q_i} M_{S,p}(f)^p \right)^{1/p} \left( \sum_i \|g\|_{W^{1,p'}(Q_i)}^{p'} \right)^{1/p'}.$$
The first term is bounded by $\|M_{S,p}(f) \|_{L^p}$.
 The second term is bounded by $\|g\|_{W^{1,p'}}=1$ since the collection $(Q_{i})_i$ forms a bounded covering. Therefore the  proposition follows.
\findem

\gb We also would like to prove a similar result as in Proposition \ref{prop:importante2} with a maximal operator
$M_{S,r}$, given by another exponent $r\leq p$. Such a result for $r\geq p$ holds combining Remark
\ref{rem:comparaison} and Proposition \ref{prop:importante2}. For $r<p$ this fact seems not to be obvious and we do
not know if it is true in a general case. That is why, we define the following assumption~:

\begin{ass} Take two exponents $\p_0,\p_1$ with $1\leq\p_0<\p_1<\infty$. Then we call $(H_{\p_0,\p_1})$ the following assumption: there exists an implicit constant such that for all functions $f\in W^{-1,\p_1}$
 \begin{equation}
\|f\|_{W^{-1,\p_1}} \lesssim \|M_{S,*,\p_0}(f)\|_{L^{\p_1}}. \tag{$H_{\p_0,\p_1}$} \label{Hp}
\end{equation}
\end{ass}

\begin{rem} \label{remf}
If $\p_0\geq \p_1$, we have seen that ($H_{\p_0,\p_1}$) is always satisfied.
 \end{rem}

\mb We finish this subsection comparing the two maximal operators $M_{S,p}$ and $M_{S,*,p}$. We have already seen in Remark \ref{rem:inegalite} that we have a pointwise inequality. We describe here a global reverse inequality.

\begin{prop} \label{prop:reverse}  Let $p\in(1,\infty)$ and $r\in[1,\infty)$. Assume that the Riemannian manifold $M$ satisfies $\mu(M)=\infty$\footnote{which is true if we assume $M$ complete
since here  the Riemannian measure is doubling.}. Then we have \be{eq4bis} \left\| M_{S,*,r}(f) \right \|_{L^p} \simeq
\left\| M_{S,r}(f) \right \|_{L^p}. \ee
The implicit constants can be chosen independently with respect to any function $f\in W^{-1,p}$.
\end{prop}

\dem  Using Remark \ref{rem:inegalite}, we just have to prove that
\be{eq103} \left\| M_{S,*,r}(f) \right\|_{L^p} \lesssim \left\| M_{S,r}(f) \right\|_{L^p}.\ee
The proof is based on a ``good lambdas'' inequality.
By classical arguments (see \cite{AM}), we just need to show the following inequality for any small enough $\gamma$ and a large enough numerical constant $K>1$~:
\begin{align}
& \mu\left(\left\{x, M_{S,*,r}(f)(x) >K\lambda,\ M_{S,r}(f)(x)\leq \gamma\lambda\right\}\right) \lesssim \hspace{2cm} & \nonumber \\
  & \hspace{6cm} \gamma \mu\left(\left\{x,M_{S,*,r}(f)(x)>\lambda\right\} \right). & \label{gl} 
\end{align}
We consider the sets
$$ B_\lambda:= \left\{ M_{S,*,r}(f)>K\lambda,\ M_{S,r}(f)\leq \gamma\lambda \right\}$$
and
$$ E_\lambda:= \left\{ M_{S,*,r}(f)>\lambda \right\}.$$
First we have $B_\lambda\subset E_\lambda$. We choose $(Q_j)_j$ a Whitney decomposition of $E_\lambda$ and write $x_j$
for a point in $4 Q_j \cap E_\lambda^c$. Let $x$ be a point in $B_\lambda \cap Q_j$. We have \be{eq105}
\inf_{f=\phi+div(\psi)} \sup_{\genfrac{}{}{0pt}{}{Q \textrm{ball}}{Q\ni x}}\  \frac{1}{\mu(Q)}
\left\||\phi|+|\psi|\right\|_{L^r(Q)} \geq K\lambda. \ee However for all ball $Q$ containing $x$ and satisfying $Q\cap
(8Q_j)^c \neq \emptyset$, the point $x_j$ belongs to $4Q$. Hence
$$ \inf_{f=\phi+div(\psi)} \sup_{\genfrac{}{}{0pt}{}{Q\ni x}{Q\cap (8Q_j)^c \neq \emptyset} } \frac{1}{\mu(4Q)}  \left\| |\phi|+|\psi|\right\|_{L^r(4Q)}\leq M_{S,*,r}(f)(x_j)\leq \lambda.$$
Therefore using $(D)$, we obtain
$$ \inf_{f=\phi+div(\psi)} \sup_{\genfrac{}{}{0pt}{}{x\ni Q}{Q\cap (8Q_j)^c \neq \emptyset} } \frac{1}{\mu(Q)}  \left\| |\phi|+|\psi|\right\|_{L^r(Q)}\lesssim \lambda.$$
Taking $K$ large enough (larger than the implicit constant in the previous inequality) yields
$$ \inf_{f=\phi_0+div(\psi_0)} \sup_{\genfrac{}{}{0pt}{}{Q \textrm{ball}}{x\in Q\subset 8 Q_j}}\  \frac{1}{\mu(Q)}  \left\|\left(\phi_0|+|\psi_0|\right) {\bf 1}_{8Q_j}\right\|_{L^{r}(Q)} \geq K\lambda.$$
Now we choose $\phi_j$ and $\psi_j$ such that~
\be{hypa} \left\| |\phi_j|+|\psi_j|\right\|_{L^r(8Q_j)} \simeq \|f\|_{W^{-1,r}(8Q_j)}. \ee
This is possible due to Proposition \ref{dualS}. We thus obtain
$$ M_{HL,r}\left(\left(|\phi_j|+|\psi_j|\right) {\bf 1}_{8Q_j}\right)(x) \geq K\lambda.$$
So we have proved that~
$$ B_\lambda \cap Q_j \subset \left\{x,\ M_{HL,r}\left(\left(\phi_j|+|\psi_j|\right) {\bf 1}_{8Q_j}\right)(x) \geq K\lambda \right\}.$$
Using the weak type $(r,r)$ of the Hardy-Littlewood maximal operator, we deduce that~
$$ \mu\left(B_\lambda \cap Q_j\right) \lesssim \frac{1}{\lambda^r} \left\||\phi_j|+|\psi_j|\right\|_{L^r(8Q_j)}^r \lesssim \frac{1}{\lambda^r} \|f\|_{W^{-1,r}(8Q_j)}^r.$$
The last inequality is due to (\ref{hypa}). Then by definition of $M_{S,r}$, we have~
\begin{align*}
 \|f\|_{W^{-1,r}(8Q_j)} & \lesssim \mu(Q_j)^{1/r} \inf_{8Q_j} M_{S,r}(f) 
 %\lesssim \mu(Q_j)^{1/r} M_{S,r}(f)(x_j) \\
 %& \lesssim \mu(Q_j)^{1/r} M_{S,*,r}(f)(x_j)
 \lesssim \gamma \mu(Q_j)^{1/r} \lambda.
\end{align*}
We conclude that
$$ \mu\left(B_\lambda \cap Q_j\right) \lesssim \gamma^r \mu(Q_j).$$
Therefore summing over $j$, the proof of (\ref{gl}) is completed. \findem

\begin{cor} \label{cor}  Let $1\leq\p_0,\p_1<\infty$ . Then under the assumptions (\ref{Hp}) and $\mu(M)=\infty$, the maximal operator $M_{S,\p_0}$ ``controls the Sobolev norm in $W^{-1,\p_1}$'': that is
\be{eq3bise} \|f\|_{W^{-1,\p_1}} \lesssim \left\| M_{S,\p_0}(f) \right \|_{\p_1}. \ee
\end{cor}

\gb It is difficult to check the assumption (\ref{Hp}), some technical details create problems. We are going to check
that the assumption (\ref{Hp}) has really a sense and is satisfied under more classical assumptions. The next
subsection is devoted to prove that (\ref{Hp}) holds under usual assumptions on the manifold $M$. This is the main
result of this section.

\subsection{Some hypotheses insuring (\ref{Hp}).}

\mb We first define some concepts to describe our main result.

\begin{df} We use the second order operator $L:=(I+\Delta)$ defined with the positive Laplace-Beltrami operator. We recall that the two operators $\Delta$ and $L$ are self-adjoint. \\
 According to \cite{AC}, we say that for $p\in(1,\infty)$ we have the non-homogeneous property (\ref{Rp}) if
 \begin{equation} \label{Rp}
 \|f\|_{W^{1,p}} \lesssim \left\| L^{1/2}(f) \right\|_{L^p} \tag{$nhR_p$}
 \end{equation}
 for all $f\in C_0^{\infty}(M)$. This is equivalent to the $L^p$ boundedness of the local Riesz transform $\nabla(\Delta+I)^{-1/2}$.\
And we have the non-homogeneous reverse property (\ref{RRp}) if
 \begin{equation} \label{RRp}
   \left\| L^{1/2}(f) \right\|_{L^p} \lesssim \|f\|_{W^{1,p}} \tag{$nhRR_p$}
 \end{equation}
  for all $f\in C_0^{\infty}(M)$.
\end{df}

\begin{df} \label{df:od} Let $p,q\in[1,\infty)$. We say that the collection $(T_t)_{t>0}=(e^{-t\Delta})_{t>0}$ or $(T_t)_{t>0}=(\sqrt{t}\nabla e^{-t\Delta})_{t>0}$  satisfy ``$(L^{p}-L^{q})$-off-diagonal estimates'', if there exists $\gamma$ such that for all balls $Q$ of radius $r_Q$, every function $f$ supported on $Q$ and all index $j\geq 0$
$$ \left( \frac{1}{\mu(2^jQ)}\int_{S_j(Q)} \left| T_{r_Q^2}(f) \right|^{q} d\mu \right)^{1/q} \lesssim  e^{-\gamma 4^j} \left( \frac{1}{\mu(Q)}\int_Q \left| f \right|^{p} d\mu \right)^{1/p}.$$
We used $S_j(Q)$ for the dyadic corona around the ball
$$S_j(Q):=\left\{y, 2^{j}\leq 1+\frac{d(y,Q)}{r_Q} <2^{j+1}\right\}.$$
These ``off-diagonal estimates'' are closely related to ``Gaffney estimates'' of the semigroup.
\end{df}

\mb We now come to our main result~:

\begin{thm} \label{poincare} Let $1<s<r'<\sigma$. Assume that the Riemannian manifold $M$ satisfies $(nhRR_{r})$ and $(nhR_{s'})$. Moreover assume that the semigroup $(e^{-t\Delta})_{t>0}$ satisfies ``$(L^{\sigma'}-L^{s'})$-off-diagonal estimates'' and that the collection $(\sqrt{t}\nabla e^{-t\Delta})_{t>0}$ satisfies ``$(L^{s'}-L^{s'})$-off-diagonal estimates''. \\
Then there is a constant $c=c(s,r,\sigma)$ such that
\be{eq3bisse} \forall f\in W^{-1,r'}, \qquad \|f\|_{W^{-1,r'}} \lesssim \left\| M_{S,*,s}(f) \right \|_{L^{r'}}. \ee
Therefore (\ref{Hp}) is satisfied for all exponents $\p_0,\p_1$ satisfying $ \p_0\geq s$ and $\p_1 =r'$.
\end{thm}

\dem  Thanks to Proposition \ref{prop:importante2}, this result is interesting only for $s<r'$, which will be assumed. \\
The proof is quite technical, we deal with the case where the manifold is of infinite measure $\mu(M)=\infty$. We explain in Remark \ref{remmf}, the modifications one has to do in the other case. \\
Take a function $f\in W^{-1,r'}$. By definition, $(nhRR_{r})$ implies that
\begin{align}
 \|f\|_{W^{-1,r'}} & := \sup_{\genfrac{}{}{0pt}{}{g\in C^\infty_0}{\|g\|_{W^{1,r}}\leq 1}}| \langle f, g\rangle| \nonumber \\
 &  = \sup_{\genfrac{}{}{0pt}{}{g\in C^\infty_0}{\|g\|_{W^{1,r}}\leq 1}} \langle L^{-1/2}f, L^{1/2} g\rangle \nonumber \\
 &  \lesssim \sup_{\|h\|_{L^{r}}\lesssim 1} \langle L^{-1/2}f, h\rangle \nonumber \\
 &  \simeq \| L^{-1/2}f\|_{L^{r'}}. \label{sob1}
\end{align}
Now we have to use a ``Fefferman-Stein'' inequality adapted to our operator $L^{-1/2}$. We use the results of \cite{BJ}. Let us first recall some notations.\\
We set
$$ M_{\sigma}(f)(x) := \sup_{Q\ni x}\ \frac{1}{\mu(Q)^{1/\sigma}} \left\| e^{-r_Q^2\Delta}(f) \right\|_{L^\sigma(Q)}$$
and
$$ M^\sharp_{s}(f)(x) := \sup_{Q\ni x}\ \frac{1}{\mu(Q)^{1/s}} \left\| f-e^{-r_Q^2\Delta}(f) \right\|_{L^s(Q)}.$$
The assumed ``$(L^{\sigma'}-L^{s'})$-off-diagonal estimates'' for $(e^{-t\Delta})_{t>0}$ gives (see \cite{BJ}, Theorem 5.11)
$$ M_{\sigma}(f) \lesssim M_{HL,s}(f).$$
Moreover from \cite{BJ}, Proposition 7.1 (which proves that the associated atomic Hardy space is included in $L^1$) and
Corollary 5.8, it comes that for all $q\in(s,\sigma)$
\begin{equation}\label{mshap} \| . \|_{L^q} \simeq \left\|{M_s}^\sharp(.) \right\|_{L^q} .
\end{equation}
We have used here that $\mu(M)=\infty$ (see Remark \ref{remmf}). \\
Thus applying (\ref{sob1}) and (\ref{mshap}) with $q=r'$, we obtain
$$ \|f\|_{W^{-1,r'}} \lesssim \left\|M_s^\sharp(L^{-1/2}f) \right\|_{L^{r'}}.$$
It remains to prove the following property \be{amont} M_s^\sharp(L^{-1/2}f) \lesssim M_{S,*,s}(f).\ee
\\
Fix an $x_0\in M$. Take a decomposition $f=\phi-div(\psi)$, with $\phi\in L^{p}$ and $\psi \in \mathcal{D}'(M)$ and a ball $Q\ni x_0$ such that
$$M_s^\sharp(L^{-1/2}f)(x_0) \leq 2 \frac{1}{\mu(Q)^{1/s}} \left\| (1-e^{-r_Q^2\Delta})L^{-1/2}f \right\|_{L^s(Q)}.$$
We can find a function $g\in C^\infty_0(Q)$ with $\|g\|_{L^{s'}}\leq 1$ such that
\begin{align*}
M_s^\sharp(L^{-1/2}f)(x_0) & \leq  \frac{4}{\mu(Q)^{1/s}} \langle (1-e^{-r_Q^2\Delta})L^{-1/2}f,g \rangle \\
 & \leq  \frac{4}{\mu(Q)^{1/s}} \langle f,L^{-1/2}(1-e^{-r_Q^2\Delta})g \rangle.
\end{align*}
Using the decomposition of $f$, we get~:
\begin{align*}
\lefteqn{M_s^\sharp(L^{-1/2}f)(x_0) \leq } & & 
 \\ & &   \frac{4}{\mu(Q)^{1/s}} \left[\langle \phi,L^{-1/2}(1-e^{r_Q^2\Delta})g \rangle + \langle \psi,\nabla L^{-1/2}(1-e^{-r_Q^2\Delta})g \rangle \right].
\end{align*}
Let us study the first term $\langle \phi,L^{-1/2}(1-e^{-r_Q^2\Delta})g \rangle$.
We follow ideas of \cite{A} (section 4, Lemma 4.4), using the following representation of the square root~:
$$ L^{-1/2}(h) = \int_0^{\infty} e^{-t} e^{-t\Delta}(h) \frac{dt}{\sqrt{t}}.$$
Now using the $(L^{s'}-L^{s'})$-``off-diagonal'' decays (implied by the
$(L^{\sigma'}-L^{s'})$- ones) of the semigroup $(e^{-t\Delta})_{t>0}$, we obtain ~:
$$ \frac{1}{\mu(2^jQ)^{1/s'}} \| L^{-1/2}(1-e^{-r_Q^2\Delta})g \|_{L^{s'}(S_j(Q))} \lesssim 2^{-j} \frac{1}{\mu(Q)^{1/s'}} \| g \|_{L^{s'}(Q)}.$$
We do not detail the proof of this claim and refer the reader to \cite{A}, Lemma 4.4 in the Euclidean case and for  $j\geq 2$. For $j\in\{0,1\}$, this is a direct consequence of the $L^{s'}$-boundedness of the semigroup.
With the normalization of $g$, we finally get
$$ \frac{1}{\mu(Q)^{1/s}}\left|\langle \phi,L^{-1/2}(1-e^{-r_Q^2\Delta})g \rangle \right| \lesssim M_{HL,s}(\phi)(x_0).$$
Similarly by the ``off-diagonal'' decays of $(\sqrt{t}\nabla e^{-t\Delta})_{t>0}$, we obtain for $j\geq 2$~:
$$ \frac{1}{\mu(2^jQ)^{1/s'}} \| \nabla L^{-1/2}(1-e^{-r_Q^2\Delta})g \|_{L^{s'}(S_j(Q))} \lesssim 2^{-j} \frac{1}{\mu(Q)^{1/s'}} \| g \|_{L^{s'}(Q)}.$$
When $j\in\{0,1\}$, this inequality is a consequence of the $L^{s'}$-boundedness of the non-homogeneous Riesz transform due to $(nhR_{s'})$. Therefore
$$ M_s^\sharp(L^{-1/2}f)(x_0) \lesssim M_{HL,s}(\phi)(x_0)+ M_{HL,s}(\psi)(x_0)$$
Taking the infimum over all the decompositions of $f$ yields (\ref{amont}) and the proof is therefore complete. \findem

\begin{rem} \label{remmf} In the case where the manifold is of finite measure, the ``Fefferman-Stein'' inequality (\ref{mshap}) has to be replaced by the following one~:
\be{msharp2} \| . \|_{L^q} \simeq \left\|M_s^\sharp(.) \right\|_{L^q} + \| . \|_{L^1}. \ee
However when $M$ is of finite measure, we have~:
$ \|L^{-1/2}(f)\|_{L^1} \lesssim \|L^{-1/2}(f) \|_{L^{s}}$. Then using the $(nhR_{s'})$ property, we deduce that
$$ \|L^{-1/2}(f)\|_{L^1} \lesssim \|f \|_{W^{-1,s}}.$$
The reverse inequality of Proposition \ref{prop:importante2} gives us
$$ \|L^{-1/2}(f)\|_{L^1} \lesssim \|f \|_{W^{-1,s}} \lesssim \|M_{S,*,s}(f)\|_{L^s}$$
which implies the desired inequality
$$ \|L^{-1/2}(f)\|_{L^1} \lesssim \|M_{S,*,s}(f)\|_{L^{r'}}$$
when $s\leq r'$.
\end{rem}

\mb We recall criterions from \cite{AC}, \cite{ACDH}, \cite{CD1}, \cite{gri} that insure our previous assumptions~:

\begin{thm}\label{rapel} Let $M$ be a complete doubling Riemannian manifold.
\begin{itemize}
 \item $(nhR_2)$ and $(nhRR_2)$ are always satisfied. 
 \\
 \item
 Assume that the heat kernel $p_t$ of the semigroup $e^{-t\Delta}$ satisfies the following pointwise estimate~
 \begin{equation}\label{due}
 p_t(x,x) \lesssim \frac{1}{\mu(B(x,t^{1/2}))} \tag{$DUE$}.
 \end{equation}
Then $(D)$ and $(DUE)$ imply the following gaussian upper-bound estimate of $p_t$
 \begin{equation}\label{ue}
 p_t(x,y) \lesssim \frac{1}{\mu(B(y,t^{1/2}))}e^{-c\frac{d^{2}(x,y)}{t}} \tag{$UE$}.
\end{equation}
Note that under $(UE)$, the collections $(e^{-t\Delta})_{t>0}$ and $(\sqrt{t}\nabla e^{-t\Delta})_{t>0}$ satisfy ``$(L^2-L^2)$ off-diagonal decays''.
Moreover we have $(L^1-L^\infty)$ ``off-diagonal'' decays of $(e^{-t\Delta})_{t>0}$.
\item It is known that the conjunction of $(D)$ and Poincar\'e inequality $(P_2)$ on $M$ is equivalent to
\begin{equation}\label{ly}
\frac{1}{\mu(B(y,t^{1/2}))}e^{-c\frac{d^{2}(x,y)}{t}}\lesssim p_t(x,y) \lesssim \frac{1}{\mu(B(y,t^{1/2}))}e^{-c\frac{d^{2}(x,y)}{t}} \tag{$LY$}.
\end{equation}
 \item Assume that  $p_t$ satisfies (\ref{due}). Then for all $p\in(1,2]$, $(nhR_p)$ and $(nhRR_{p'})$ hold. Moreover, $p_t$ satisfies a gaussian upper-bound estimate, \item Under Poincar\'e inequality $(P_2)$, the property $(nhR_p)$ for all $p\in(2,p_0)$ is equivalent to the boundedness
\begin{equation}\label{gp}
 \left\|\nabla e^{-t\Delta} \right\|_{L^{p_0} \to L^{p_0}} \lesssim \frac{1}{\sqrt{t}}\tag{$G_{p_0}$}.
 \end{equation}
Moreover $(G_{p_0})$ implies for $p\in[2,p_0)$ the $(L^p-L^p)$ ``off-diagonal'' decays of $(\sqrt{t}\nabla e^{-t\Delta})_{t>0}$.
 \item Under Poincar\'e inequality $(P_{p_0})$ for $p_0\in (1,2]$, $(nhRR_p)$ holds for all $p\in (p_0,2]$.
\end{itemize}
\end{thm}

\begin{rem} All these results are proved in their homogeneous version, with homogeneous properties $(R_p)$ and $(RR_p)$. It is
essentially based on the well-known Calder\'on-Zygmund decomposition for Sobolev functions. This tool was extended for
non-homogeneous Sobolev spaces (see \cite{Nadine}). Thus by exactly the same proof, we can obtain an analogous
non-homogeneous version  and then prove Theorem \ref{rapel}.
\end{rem}

\mb From Theorem \ref{poincare}, Theorem \ref{rapel}, Remark \ref{rem:comparaison} and the self-improvement of Poincar\'e inequality (proved in \cite{KZ}) we get~:

\begin{cor} \label{corpoi} Let $M$ be a non-compact Riemannian manifold satisfying the doubling property.
If Poincar\'e inequality $(P_{p_0})$ holds for some $p_0\in (1,2)$, then (\ref{Hp}) is verified for all $\p_0,\p_1$
satisfying
$$ \p_0\geq 2 \qquad \textrm{and}  \qquad \p_1 \leq p_0'.$$
\end{cor}

\begin{cor} \label{corRn} In the Euclidean case $M=\R^n$, for all $\p_0,\p_1\in(1,\infty)$, the assumption (\ref{Hp}) holds. More generally, on any Riemannian manifold satisfying $(D)$ and $(P_1)$, (\ref{Hp}) holds for all $\p_0,\p_1\in(1,\infty)$.
\end{cor}

\mb We begin to understand the link between Sobolev norms and the Lebesgue norms of our maximal operators. This technical result will be useful in Section \ref{sec:interpolation} to develop new results for the interpolation of Sobolev spaces.

\section{Interpolation of Sobolev spaces.} \label{sec:interpolation}

In this section, we look for a real interpolation result for the scale of Sobolev spaces $(W^{1,p})_{p\in(1,\infty)}$.
We refer the reader to the work of N.~Badr (see \cite{Nadine,Nadine1}) for first results. This work is based on a well-known Calder\'on-Zygmund decomposition for Sobolev functions, initialy explained by P. Auscher in \cite{A1}. We refer the reader to \cite{A1} for the first use of this one. Many applications follow from this decomposition and there are many versions (for example there is several improvements with weights in \cite{AM} and \cite{Nadine1}). This very useful tool works under the assumption of Poincar\'e inequality.

\gb This section is devoted to the description of interpolation results for Sobolev spaces using the results of Section
\ref{sec:maximal}.

\gb We recall the important assumption~:

\begin{ass} Take two exponents $1\leq\p_0\leq \p_1<\infty$. We call $(H_{\p_0,\p_1})$ the following assumption~:
 \begin{equation}
\|f\|_{W^{-1,\p_1}} \lesssim \|M_{S,*,\p_0}(f)\|_{L^{\p_1}}. \tag{$H_{\p_0,\p_1}$} \label{Hp}
\end{equation}
\end{ass}

\begin{df} \label{ensembleI} For $M$ a Riemannian manifold, we denote by ${\mathcal I}_M$ the following set
:
$$ {\mathcal I}_M:=\left\{ (\p_0,\p_1)\in(1,\infty)^2,\ \p_0\leq \p_1,\ (\textrm{\ref{Hp})} {\textrm{ holds}} \right\}.$$
\end{df}

\mb Here is the main result of this subsection~:

\begin{thm} \label{thm:inter}Let $M$ be a Riemannian manifold satisfying the doubling property $(D)$. Then
the scale $(W^{1,p})_{p\in(1,\infty]}$ is an interpolation scale for the real interpolation related to ${\mathcal I}_M$. That is for all $p_0,p_1\in(1,\infty]$ (with $p_0\leq p_1$) and $\theta\in(0,1)$ such that
$$ \frac{1}{p_\theta}:=\frac{1-\theta}{p_0}+\frac{\theta}{p_1}$$
and satisfying $(p_1',p_\theta')\in {\mathcal I}_M$, we have
$$\left(W^{1,p_0},W^{1,p_1}\right)_{\theta,p_\theta}=W^{1,p_\theta}.$$
\end{thm}

\dem We set $E:=\left(W^{1,p_0},W^{1,p_1}\right)_{\theta,p_\theta}$. We have to prove the equivalence of norms~:
\be{equiv} \|\ \|_{E} \simeq \|\ \|_{W^{1,p_\theta}}. \ee
>From the interpolation theory on Lebesgue spaces, it is obvious that
$$ E \hookrightarrow W^{1,p_\theta}.$$
We just have to prove the reverse embedding. We will use the maximal operator $M_{S,*,p_1'}$.
Let $q\in[1,p_1]$. We claim that there is a constant $c=c_q$ such that~
\be{claim1} \left\| M_{S,*,p_1'}(h) \right\|_{L^{q',\infty}} \leq c_q \|h\|_{(W^{1,q})^*}. \ee
This fact comes from several properties~: $(W^{1,q})^*=W^{-1,q'}$ by definition, from (\ref{eq1}) and finally from $M_{S,*,p_1'}\leq M_{S,*,q'}$. \\
We refer the reader to \cite{KP} for the proof of
$$ E^*=\left[(W^{1,p_0},W^{1,p_1})_{\theta,p_\theta} \right]^* = (W^{-1,p_0'},W^{-1,p_1'})_{\theta,p_{\theta}'} $$
with the concept of ``doolittle couple''. \\
Then by interpolation on the weak Lebesgue spaces, we obtain~:
$$ \left\| M_{S,*,p_1'}(h) \right\|_{L^{p_\theta'}} \lesssim \|h\|_{E^*},$$
which according to the assumption (\ref{Hp}) (for $\p_0=p_1'$ and $\p_1=p_\theta'$) yields
$$ \left\| h \right\|_{W^{-1,p_\theta'}} \lesssim \|h\|_{E^*}.$$
Thus $E\hookrightarrow W^{1,p_\theta}$ and $E^* \hookrightarrow W^{-1,p_\theta'}=(W^{1,p_\theta})^*$. Using Hahn-Banach Theorem, we deduce that $E=W^{1,p_\theta}$ with equivalent norms.
\findem

\mb To regain results of the same kind as in \cite{Nadine}, where the author assumes Poincar\'e inequality, we describe
the following corollary~:

\begin{cor} \label{cor2} Assume that $M$ satisfies $(D)$ and admits a Pincar\'e inequality $(P_r)$ for an $r\in(1,2)$.
Then for all $p_0\in(1,2)$ and $\theta\in(0,1)$ such that
$$ \frac{1}{p_\theta}:=\frac{1-\theta}{p_0}+\frac{\theta}{2} <\frac{1}{r},$$
we have
$$\left(W^{1,p_0},W^{1,2}\right)_{\theta,p_\theta}=W^{1,p_\theta}.$$
\end{cor}

\dem We set $\p_1=p_\theta'$ and $\p_0=2$. Thanks to Theorem \ref{thm:inter}, we just have to check that $(\p_0,\p_1)\in {\mathcal I}_M$.
This is a direct consequence of Corollary \ref{corpoi}. \findem

\begin{rem} In the previous corollary, Poincar\'e inequality $(P_r)$ could be replaced by the weaker non-homogeneous variant 
\begin{equation*}\tag{$\tilde{P}_{r}$}
\left(\aver{Q}|f-f_{Q}|^{r} d\mu\right)^{1/r} \leq C r_Q \left(\aver{Q}(|f|^r+|\nabla f|^{r}d\mu\right)^{1/r}.
\end{equation*}
As ($\tilde{P}_{r}$) is sufficient to obtain the ($nhRR_r$) property.
Moreover Assumption $(nhRR_{r})$ is sufficient for the previous corollary.
\end{rem}

\begin{rem} In Corollary \ref{cor2}, we can chose $p_1\leq 2$ (and not necessary equal to $2$). Then under $(DUE)$ and $(nhRR_r)$, we get the corresponding interpolation result.
\end{rem}

\mb Let us compare these results with \cite{Nadine}. Note first that the results --even the proofs-- of \cite{Nadine} in the non-homogeneous case still hold with this variant of Poincar\'e inequality. In \cite{Nadine}, the author just requires the condition $p_\theta> r$ to obtain the interpolation result under
{\it{local}} doubling property and {\it{local}} Poincar\'e inequalities. The main tool (the ``well-known'' Calder\`on-Zygmund decomposition for Sobolev functions) of \cite{Nadine} permits to interpolate any Sobolev spaces (not only with $W^{1,2}$ or $W^{1,p_1}$ with $p_1\leq 2$) under Poincar\'e inequality $(P_r)$. \\ 
The use of the exponent $2$ is the most important in the litterature and that is why we mainly deal with it. In the case $p_1\leq 2$, our assumption $(nhRR_{r})$ is weaker than the corresponding Poincar\'e inequality $(P_r)$. Consequently we regain the results of N. Badr (\cite{Nadine}). \\
However in the case where $p_1>2$, we can not recover her results as we require an extra assumption~: the Riesz inequality. Our assumptions and the ones of \cite{Nadine} are not comparable when $p_1>2$.
Which is interesting is that even in this case, we succeed to interpolate Sobolev spaces without assuming Poincar\'e inequalities. \\
An interesting question still stays open~: we have weaken the assumption of Poincar\'e inequality, however we do not know which assumptions should be sufficient and necessary to prove an interpolation result. In the case $p_0,p_1\leq 2$ our assumption $(nhRR_r)$ seems to be the well-adapted assumption ...

\gb
To finish, we refer the reader to an other work (joined with N.~Badr, see \cite{HSNadine}), where we develop a new theory for abstract Hardy-Sobolev spaces. Using these maximal inequalities, we prove some results for interpolation between Hardy-Sobolev spaces and Sobolev spaces. In this application, the arguments based on the well-known Calder\`on-Zygmund decomposition do not work and these new maximal inequalities play a crucial role. 

%\section*{Acknowledgements} I am sincerely indebted to Nadine Badr, who has carefully read and who proposed to me valuable suggestions and corrections.

\bibliographystyle{plain}
\bibliography{biblio}

\begin{thebibliography}{10}

\bibitem{ambrosio1}
L.~Ambrosio, Jr.~M. Miranda, and D.~Pallara.
\newblock Special functions of bounded variation in doubling metric measure
  spaces.
\newblock {\em Calculus of variations~: topics from the mathematical heritage
  of E. De Giorgi, Quad. Mat., Dept. Math, Seconda Univ. Napoli, Caserta
  {\textbf 14}}, pages 1--45, 2004.

\bibitem{A1}
P.~Auscher.
\newblock On ${L}^p$ estimates for square roots of second order elliptic
  operators on ${\R}^{n}$.
\newblock {\em Publ. Mat. \textbf{48}}, pages 159--186, 2004.

\bibitem{A}
P.~Auscher.
\newblock On necessary and sufficient conditions for ${L}^p$ estimates of
  {R}iesz transforms associated to elliptic operators on ${\R}^n$ and related
  estimates.
\newblock {\em Memoirs of Amer. Math. Soc. \textbf{186} no.871}, 2007.

\bibitem{AC}
P.~Auscher and T.~Coulhon.
\newblock Riesz transform on manifolds and {P}oincar{\'e} inequalities.
\newblock {\em Ann. Sc. Nor. Sup. Pisa (\textbf{5}), IV, 3}, pages 531--555,
  2005.

\bibitem{ACDH}
P.~Auscher, T.~Coulhon, X.T. Duong, and S.~Hofmann.
\newblock Riesz transform on manifolds and heat kernel regularity.
\newblock {\em Ann. Sci. Ecole Norm. Sup. {\textbf{37}}}, pages 911--957, 2004.

\bibitem{AM}
P.~Auscher and J.M. Martell.
\newblock Weighted norm inequalities, off-diagonal estimates and elliptic
  operators, {Part I} : {G}eneral operator theory and weights.
\newblock {\em Adv. in Math. no {\textbf{212}}}, pages 225--276, 2007.

\bibitem{ART}
P.~Auscher, E.~Russ, and P.~Tchamitchian.
\newblock Hardy-sobolev spaces on strongly lipschitz domains of ${\R}^n$.
\newblock {\em J. Func. Anal. \textbf{218}}, pages 54--109, 2005.

\bibitem{Nadine}
N.~Badr.
\newblock Real interpolation of {S}obolev spaces.
\newblock {\em Mathematica Scandinavica}, 2008.

\bibitem{Nadine1}
N.~Badr.
\newblock Real interpolation of {S}obolev spaces associated to a weight.
\newblock {\em submitted}, 2008.

\bibitem{HSNadine}
N.~Badr and F.~Bernicot.
\newblock Abstract hardy-sobolev spaces and interpolation.
\newblock {\em submitted}, 2008.

\bibitem{B2}
F.~Bernicot.
\newblock Use of abstract {H}ardy spaces, real interpolation and applications
  to bilinear operators.
\newblock {\em submitted}, page available at http://fr.arxiv.org/abs/0809.4110,
  2008.

\bibitem{BJ}
F.~Bernicot and J.~Zhao.
\newblock Abstract {H}ardy spaces.
\newblock {\em J. Funct. Anal. \textbf{255} no. 7}, pages 1761--1796, 2008.

\bibitem{coifman2}
R.~Coifman and G.~Weiss.
\newblock {\em Analyse harmonique sur certains espaces homog\`{e}nes}.
\newblock Lecture notes in Math., Springer, 1971.

\bibitem{CD1}
T.~Coulhon and X.T. Duong.
\newblock Riesz transforms for $1\leq p\leq 2$.
\newblock {\em Trans. Amer. Math. Soc. \textbf{351} no.2}, pages 1151--1169,
  1999.

\bibitem{gri}
A.~Grigor'yan.
\newblock Gaussian upper bounds for the heat kernel on arbitrary manifolds.
\newblock {\em J. Diff. Geom. \textbf{45}}, pages 33--52, 1997.

\bibitem{hajlasz4}
P.~Hajlasz and P.~Koskela.
\newblock {S}obolev met {P}oincar\'{e}.
\newblock {\em Mem. Amer. Math. Soc. \textbf{145} no. 688}, pages 1--101, 2000.

\bibitem{KP}
S.~Kaijser and J.~W. Pelletier.
\newblock Interpolation {F}unctors and {D}uality.
\newblock {\em Lectures Notes Math. no.{\textbf{1208}}}, 1980.

\bibitem{KZ}
S.~Keith and X.~Zhong.
\newblock The {P}oincar{\'e} inequality is an open ended condition.
\newblock {\em Ann. of Maths. \textbf{167} no.2}, pages 575--599, 2008.

\end{thebibliography}

\end{document}